\documentclass[letterpaper,12pt]{article}
\usepackage{amssymb}
\usepackage{amsmath}
\usepackage{amscd}
\usepackage[all]{xy}

\usepackage{amsthm}

\usepackage{centernot}

\usepackage{bm}
\usepackage{setspace}

\newtheorem{theo}{Theorem}[section]
\newtheorem{prop}[theo]{Proposition}
\newtheorem{lem}[theo]{Lemma}
\newtheorem{cor}[theo]{Corollary}
\newtheorem{defi}[theo]{Definition}

\def \Br {{\rm{Br}}}

\def \Pic {{\rm {Pic}}}

\def \Gal {{\rm{Gal}}}

\def \Im {{\rm {Im\,}}}

\def \A{{\bm A}}

\def \P{{\mathbb P}}
\def \Spec {{\rm{Spec\,}}}

\def \Hom {{\rm {Hom}}}

\def \Pic {{\rm {Pic}}}

\def\ov{\overline}

\def \Z {{\mathbb Z}}
\def \Q {{\mathbb Q}}

\def\O{{\cal O}}

\def\inv{{\rm inv}}

\def\O{{\cal O}}

\newcommand{\bthe}{\begin{theo}}
\newcommand{\ble}{\begin{lem}}
\newcommand{\bpr}{\begin{prop}}
\newcommand{\bco}{\begin{cor}}
\newcommand{\bde}{\begin{defi}}
\newcommand{\ethe}{\end{theo}}
\newcommand{\ele}{\end{lem}}
\newcommand{\epr}{\end{prop}}
\newcommand{\eco}{\end{cor}}
\newcommand{\ede}{\end{defi}}

\def\beq{\begin{equation} \label}

\DeclareFontFamily{U}{wncy}{}
    \DeclareFontShape{U}{wncy}{m}{n}{<->wncyr10}{}
    \DeclareSymbolFont{mcy}{U}{wncy}{m}{n}
    \DeclareMathSymbol{\Sh}{\mathord}{mcy}{"58}

\usepackage{tikz-cd}

\usepackage[greek,english]{babel}
\usepackage[iso-8859-7]{inputenc}

\title{The Brauer-Manin obstruction for zero-cycles on $K3$ surfaces}
\author{Evis Ieronymou}
\date{}

\begin{document}
\baselineskip=15pt

\maketitle

\begin{abstract}
We study local-global principles for zero-cycles on $K3$ surfaces defined over number fields. We follow an idea of Liang to use the trivial fibration over the projective line.

\end{abstract}

\section{Introduction}

Let $X$ be a smooth, projective, geometrically irreducible variety over a number field $k$. The theory of the Brauer-Manin obstruction to the Hasse principle, and its variants, has developed a lot since its introduction by Manin
in \cite{Man} (see eg \cite{Sk01} and references therein).
Primarily one is interested in rational points but there is an analogous story for zero-cycles, where one basically replaces rational points with zero-cycles.
 This note is concerned with local-global principles for zero-cycles on $K3$ surfaces over number fields.

 First some general context. Denote by $CH_0(X)$ the group of zero-cycles modulo rational equivalence on $X$ and by $A_0(X)$ the kernel of the degree map
$CH_0(X)\to \Z$. We have a pairing
 $$
  \Br(X) \times \prod_{v\in \Omega}CH_0(X_{k_v}) \to \Q/\Z$$
  and the image of $CH_0(X)\to \prod_{v\in \Omega}CH_0(X_{k_v})$ is orthogonal to $\Br(X)$.

  The conjecture that the Brauer-Manin obstruction to the existence of a zero-cycle of degree one is the only one is still open.
    This conjecture was put forward in various forms and for various classes of smooth, projective, geometrically irreducible varieties by Colliot-Th\'el\`ene, Sansuc, Kato and Saito (cf. \cite{CTS81}, \cite{KS86}, \cite{CT95}).
    We now have the following for arbitrary smooth, projective, geometrically connected varieties over a number field (we refer the reader to \cite{CT95} and \cite[\textsection 0,1]{Wit12} for more information)

CONJECTURE $(E_0)$

The sequence
$$
\varprojlim_n A_0(X)/n  \to   \prod_{v\in \Omega}\varprojlim_n  A_0(X_{k_v})/n \to \Hom(\Br(X), \Q/\Z)
$$
 is exact.

 CONJECTURE $(E_1)$

If there is a family of local zero-cycles $\{z_v\}_{v \in \Omega}$ of degree one orthogonal to $\Br(X)$ then there exists a zero-cycle of degree one on $X$.

CONJECTURE $(E)$

 The sequence
$$
\widehat{CH_0(X)}  \to \widehat {CH_{0,\A}(X)} \to \Hom(\Br(X), \Q/\Z)
$$
 is exact (the notation will be explained below).

 This formulation is due to \cite{vHam03}, see also \cite{Wit12}.
 Note that CONJECTURE $(E)$ implies the other two (\cite[Rem. 1.1, (ii)-(iii)]{Wit12})

Our results are in the spirit of \cite{Lia}: we assume results about rational points for all finite extensions and deduce results about zero-cycles. We remind the reader that the phrase ``the Brauer-Manin obstruction to the Hasse principle on $X$ is the only one'' is the statement ``$X(\A_k)^{\Br(X)}\neq\emptyset \Rightarrow X(k)\neq\emptyset$''.
We can now recall \cite[Thm. 3.2.1]{Lia} which states the following:

Let $Y$ be a proper smooth variety defined over a number field $k$ such
that $NS(Y_{\ov k})$ is torsion-free and $H^1(Y_{\ov k},\mathcal O_{Y_{\ov k}}) = H^2(Y_{\ov k}, \O_{Y_{\ov k}}) = 0$.  Suppose that for all finite extensions $F/k$ the Brauer-Manin obstruction to the Hasse principle on $Y_F$ is the only one.
Then the Brauer-Manin obstruction is the only obstruction to the Hasse principle for $0$-cycles of degree $1$ on $Y_k$. There is also a similar statement concerning weak approximation.

Note that a rationally connected variety satisfies the conditions above, but for a $K3$ surface the assumption that $H^2(Y_{\ov k}, \O_{Y_{\ov k}}) = 0$ fails.
In this note we show that we can circumvent this crucial assumption by using results of Skorobogatov and Zarhin \cite{SZ08} and Orr and Skorobogatov \cite{OS}, and therefore extend some of the main results of \cite{Lia} to the case of $K3$
surfaces. More precisely we have:

 \bthe\label{mainthm}
Let $X$ be a $K3$ surface over a number field $k$ and $d\in\Z$. Suppose that for all finite extensions $F/k$ the Brauer-Manin obstruction to the Hasse principle on $X_F$ is the only one.
Then the Brauer-Manin obstruction to the existence of a zero-cycle of degree $d$ on $X$ is the only one: if there is a family of local zero-cycles $\{z_v\}_{v \in \Omega}$ of degree $d$ orthogonal to $\Br(X)$ then there exists a
zero-cycle of degree $d$ on $X$. \ethe

We also show something about ''weak approximation''

\bthe\label{WA}
Let $X$ be a $K3$ surface over a number field $k$. Suppose that for all finite extensions $F/k$ we have that $X(F)$ is dense in $X(\A_F)^{\Br(X_F)}$.
Let $\{z_v\}_{v \in \Omega}$ be a family of local zero-cycles  of constant degree orthogonal to $\Br(X)$.
Then for any $n\geq 1$ there exists a zero-cycle $b$ on $X$ of the same degree such that
$b=z_v$ in $CH_0(X_{k_v})/n $ for all $v\in \Omega$.

\ethe

{\it Remarks}
\begin{enumerate}
\item

The Brauer-Manin obstruction does not explain all the failures of the Hasse principle ( see \cite{Sk99} or \cite{Poo10}, where in the latter the insufficiency of an even more refined obstruction is shown). However this is
conjectured to be the case for various classes of varieties, and what is relevant to this note is that one of these classes is the class of $K3$ surfaces (cf. \cite{SZ08}, \cite{Sk09}). Note that although there is some scarce evidence for the aforementioned conjecture for rational points on $K3$ surfaces, there is no evidence whatsoever for the analogous conjecture for zero-cycles on $K3$ surfaces. By relating the two conjectures our results can be seen as providing some such evidence.

\item
The reason that Theorem \ref{mainthm} is not a Corollary of Theorem \ref{WA} is that in principle the statement ``$X(\A_k)^{\Br(X)}\neq\emptyset \Rightarrow X(k)\neq\emptyset$''  is weaker than the statement ``$X(k)$ is
dense in $X(\A_k)^{\Br(X)}$''

\item

In order to make a precise link between CONJECTURE (E) and the conclusion of Theorem
\ref{WA} we would need more information about divisibility properties of the groups $A_0(X_{k_v})$.
For example CONJECTURE (E) would follow from the conclusion of Theorem
\ref{WA} if we knew that there is a constant $C$ such that for all $v$, we have
that $C$ annihilates $A_0(X_{k_v})$ modulo its maximal divisible
subgroup. Note that this is true if $X$ is rationally connected but it is
not known in the case of $K3$ surfaces (cf. Proposition \ref{Saito} and the references in its proof).

\end{enumerate}

The plan of the proof is as follows. We employ the idea of Liang \cite[proof of Thm. 3.2.1]{Lia} to use the trivial fibration over the projective line, and hence utilize the fibration method.
 To this effect, it will be most convenient for us to quote and use some results of Harpaz and Wittenberg on the fibration method from \cite{HW}.
As already mentioned above the difficulty is that we have to somehow circumvent the crucial assumption that $H^2(\ov X,\mathcal O_{\ov X})=0$, which fails when $X$ is a $K3$ surface.
In order to do so we use some boundedness results of Skorobogatov and Zarhin \cite{SZ08} and Orr and Skorobogatov \cite{OS} on Brauer groups of $K3$ surfaces over number fields.
The organisation of the note is as follows. In \textsection 2 we set notation and conventions, in \textsection 3 we prove some auxiliary lemmata and in \textsection 4 we give the proofs of the two theorems.

\section{Notation and conventions}

All cohomology groups are \'etale (or Galois) cohomology groups. When $k$ is a field we denote by $\ov k$ a separable closure of
$k$. Let $X/k$ be a smooth, projective and geometrically irreducible variety.

We set $\Br(X)=H^2(X,\mathbb G_m)$, and $\Br_0(X)=\Im(\Br(k)\to\Br(X))$.
For $L/k$ a field extension we denote $X_L$ the variety $X \times_{\Spec(k)} \Spec(L)$ over $L$.
When $L =\ov k$ we simply write $\ov X$ . We denote by $X(L)$ the set of $L$-valued points of X, that is $X(L) =\Hom_{\Spec(k)}(\Spec(L),X)$.

When $k$ is a number field, we denote by $\Omega_k$ the set of places of $k$, and by $\A_k$ the ring of adeles. We suppress the dependence on $k$ when there is no risk of confusion. By $\Omega_f$ (resp. $\Omega_{\infty}$) we
denote the subset of $\Omega$ consisting of finite (resp. infinite) places. For $v\in \Omega$, $k_v$ is the completion of $k$ at the place $v$. For a finite subset $S\subset \Omega$ we denote by $\mathcal O_S$ the ring of
$S$-integers of $k$ (elements of $k$ that are integral at all the finite places not in $S$).

We denote by $Z_0(X)$ the group of zero-cycles of $X$, and by $CH_0(X)$ its quotient by the subgroup of cycles rationally equivalent to zero. We denote by $A_0(X)$ the kernel of the degree map
$CH_0(X)\to \Z$. A reduced zero-cycle is one all of whose coefficients belong to $\{0,1\}$.

We now follow the notations of \cite{HW}. We reproduce it here for the convenience of the reader. For us all varieties are proper, which allows for some simplifications -and so things here are actually presented in a less general
form than in \cite{HW} (see also the remark at the end of this section).

We denote by $Z_{0,\A}(X)$ the product $\prod_{v\in \Omega}Z_0(X_{k_v})$.
For the definition of $CH_{0,\A}(X)$ we modify the components at the infinite places.  We define $CH_0'(X_{k_v})$ to be the cokernel of the norm map $CH_0(X_{\ov{k_v}})\to CH_0(X_{k_v})$ when $u$ is infinite and $CH_0(X_{k_v})$
when $v$ is finite.  We denote by $CH_{0,\A}(X)$ the product $\prod_{v\in \Omega}CH'_0(X_{k_v})$.
For an abelian group $M$ we denote $\hat{M}$ the inverse limit $\varprojlim_{n\geq1}M/nM$.

We have a pairing
$$
\Br(X)\times \widehat{CH_{0,\A}(X)}\to \Q/\Z
$$

When we have a morphism $f:X\to C$ to a smooth projective curve over $k$, we denote $Z^{\text{eff,red}}_{0}(X)$
the set of effective zero-cycles $z$ on $X$ such that $f_*z$ is a reduced divisor on $C$ (by a slight abuse we suppress $f$ from the notation).
If $k$ is a number field we let $Z^{\text{eff,red}}_{0,\A}(X)=\prod_{v\in \Omega} Z^{\text{eff,red}}_{0,\A}(X_{k_v})$.
For $y\in CH_0(C)$, let $y'$ be its image in $CH_{0,\A}(C)$ and denote by $Z_{0,\A}^{\text{eff,red},y}$ the inverse image of $y'$ by the push-forward map $Z^{\text{eff,red}}_{0,\A}(X)\to CH_{0,\A}(C)$

We denote $Sym_{X/k}$ the disjoint union of the symmetric products $Sym^d_{X/k}$ for $d\geq1$.

We say that a subset $H$ of an
irreducible variety $X$ is a Hilbert subset if there exist a dense open subset
$X^0\subset X$, an integer $n\geq 1$ and irreducible finite \'etale $X^0$-schemes $W_1,\cdots,W_n$
such that $H$ is the set of points of $X^0$ above which the fiber of $W_i$ is irreducible
for all $1\leq i \leq n$. Note that by this definition $H$ is not a subset of $X(k)$.

{\it Remark}.
In this note we will only consider the trivial fibration over $\P_k^1$.
We will thus use only a small fraction of the results from \cite{HW}, which allows for some simplifications in the notation above.
For example $C$ is always $\P^1_k$ and we can take $y\in \Z$ as we identify $\Pic(\P^1)$ with $\Z$ via the degree map.
Moreover we can always take the open subsets involved in the results we will use from \cite{HW} to be the whole spaces (which are proper), which enables us to avoid the slightly more convoluted definition for $Z_{0,\A}(X)$ that
appears in op. cit.

\section{Preparatory results}

We start with a trivial lemma.
\ble\label{lempr}
Let $h:Y\to V$ be a morphism of varieties over $k$, which admits a section and let $d\in\Z$. If the Brauer-Manin obstruction to the existence of a zero-cycle of degree $d$ on $Y$ is the only one, then the same is true for $V$.
\ele
\begin{proof}
Let $f:V\to X$ be a morphism of varieties over a field $K$, and denote $f^*:\Br(X)\to\Br(V)$, and $f_*:Z_0(V)\to Z_0(X)$
the corresponding induced maps. From the functoriality of the pairings we have
$$
<\mathcal A, f_* z>=<f^*(\mathcal A),z>\in \Br(K)
$$
for $z\in Z_0(V)$ and $\mathcal A \in \Br(X)$ (cf. \cite[\textsection 3]{CTSD}).

Let $s:V\to Y$ be a section of $h$. If we have a family of local zero cycles of degree $d$ on $V$ which is orthogonal to $\Br(V)$ we can push it forward via $s$ and get a family of local zero cycles of degree $d$ on $Y$ which is
orthogonal to $\Br(Y)$. By assumption there exists a zero cycle of degree $d$ on $Y$. By pushing it forward via $h$ we get a zero cycle of degree $d$ on $V$.

\end{proof}

We will need the following easy lemma.
\ble\label{lemhs}
Let $F/k$ be a finite field extension inside $\ov k$. There exists a Hilbert subset $H$ of $\P^1_k$ with the property that for any closed point $P\in H$  we have that $F$ and $k(P)$ are linearly disjoint over $k$.
%and any embedding of $f:k(P)\to \ov k$is not contained in $f(k(P))$.
\ele
\begin{proof}
Consider the natural projection $p:\mathbb{A}^1_{F}\to \mathbb{A}^1_k$. This is a finite \'etale morphism and hence defines a Hilbert subset $H$ of $\P^1_k$. If $P$ is a closed point in $H$ then by definition $k(P)\otimes_k F$ is
irreducible and hence a field. This implies that for any $k$-linear embedding $f:k(P)\to \ov k$ we have that $f(k(P))$ and $F$ are linearly disjoint over $k$ %we have that $f(k(P))\cap L=k$

\end{proof}
We will also need the following lemma whose proof is an adaptation of the proof of \cite[Prop. 3.1.1]{Lia} to our case.

\ble\label{theL}
Let $X$ be a $K3$ surface over a number field $k$. There exists a finite Galois extension $L/k$ such that $\Gamma_L$ acts trivially on $\Pic (\ov X)$ and for any Galois extension $K/k$ that contains $L$ we have an exact
sequence
$$
0\to H^1(\Gal(K/k), \Pic \ov X)\to \Br(X)/\Br_0(X)\to \Br( \ov X)^{\Gamma_k}\to H^2(\Gal(K/k), \Pic(\ov X))
$$
\ele
\begin{proof}
From the Hochschild-Serre spectral sequence

$$
H^p(k,H^q(\ov X, \mathbb G_m))\implies  H^{p+q}(X,\mathbb G_m)
$$
and using the fact that $H^3(k,\ov k)=0$ we get
$$
0\to H^1(k, \Pic \ov X)\to \Br(X)/\Br_0(X)\to \Br( \ov X)^{\Gamma_k}\to H^2(k,\Pic( \ov X))
$$
By \cite[Thm. 1.2]{SZ08}  we have that $\Br( \ov X)^{\Gamma_k}$ is finite and hence so is its image in $H^2(k,\Pic( \ov X))$. As $H^2(k,\Pic( \ov X))$ is the direct limit of the inflation maps from finite quotients we can choose a
finite Galois extension $L_1/k$ so that
the image of $\Br( \ov X)^{\Gamma_k}$ is contained in the image of $H^2(\Gal(L_1/k),\Pic( \ov X)^{\Gamma_{L_1}})$ in $H^2(k,\Pic( \ov X))$.

As $\Pic(\ov X)$ is a finitely generated, free abelian group, we can choose a finite extension $L_2/k$ so that $\Gamma_{L_2}$ acts trivially on $\Pic(\ov X)$. Let $L$ be the Galois closure of the composite of $L_1$ and $L_2$. The
inflation-restriction sequence tells us that
 $H^1(\Gal(L/k), \Pic(\ov X))\to H^1(k, \Pic \ov X)$ is an isomorphism and $H^2(\Gal(L/k), \Pic(\ov X))\to H^2(k,\Pic( \ov X))$ is an injection. We can therefore rewrite the exact sequence above as
 $$
0\to H^1(\Gal(L/k), \Pic(\ov X))\to \Br(X)/\Br(k)\to \Br( \ov X)^{\Gamma_k}\to H^2(\Gal(L/k), \Pic(\ov X))
$$
It is clear that we can replace $L$ in the above exact sequence with $K$ where $K$ is any Galois extension $K/k$ that contains $L$.

\end{proof}

 We will moreover use the following proposition, which follows from results of Kato, Saito and Sato.

\bpr\label{Saito}
Let $V$ be a $K3$ surface over $k_v$, and let $p$ be the characteristic of the residue field. Assume that $V$ has a smooth projective model over $\O_v$. Then $A_0(V)$ is $\ell$-divisible for all primes $\ell$ with $\ell\neq p$.
\epr

\begin{proof}
This follows immediately from \cite[Thm. 0.3, Cor. 0.10]{SS10} together with  \cite[Thm. 1]{KS83}.

\end{proof}
\section{Proofs of the main results}

{\it Proof of Theorem \ref{mainthm}: } Let $Y=X\times \P^1$ and let $p:Y\to X$ be the projection to the first factor. By lemma \ref{lempr} it suffices to prove that the Brauer-Manin obstruction to the existence of a zero-cycle of
degree $d$ on $Y$ is the only one.

We start with some standard arguments (cf. \cite[\textsection 3, \textsection
5]{CTSD}). We can find a finite set $S\subset \Omega$ containing the infinite
places so that $Y$ has a smooth projective model $\mathcal Y$ over $\mathcal
O_S$, and $Y(k_v)\neq \emptyset$ for any $v\notin S$. Since $|\Br(X)/\Br_0(X)|$
is finite by \cite[Thm. 1.2]{SZ08} and $p^*:\Br(X)\to\Br(X\times \P^1)$ is an
isomorphism we have that $|\Br(Y)/\Br_0(Y)|$ is finite. Let $\alpha_i\in
\Br(Y)$, $1\leq i \leq n$, be elements whose classes generate
$\Br(Y)/\Br_0(Y)$. By enlarging $S$ if necessary, we can furthermore assume
that each $a_i$ extends to an element of $\Br(\mathcal Y)$. Hence for any
$1\leq i \leq n$, $v \notin S$ and $b\in Z_0(Y_{k_v})$ we have that
$\inv_v(<a_i,b>)=0$. Let $B$ be the span of
the $a_i$ in $\Br(Y)$.

We now work with the trivial fibration $f:Y\to \P^1$. Let $z_{\A}=\{z_v\}_{v \in \Omega}$ be a family of local zero-cycles of
degree $d$ on $Y$, which is orthogonal to $B$. By \cite[Prop. 7.5]{HW}
we can find $z'_{\A}\in Z^{\text{eff,red},\delta}_{0,\A}(Y)$, $z\in Z_0(Y)$
and $\xi_{\A}\in CH_{0,\A}(Y)$ such that the equality
$z_{\A}=z'_{\A}+z+\xi_{\A}$ holds in $ CH_{0,\A}(Y)$, the $v$-adic component
of $\xi_{\A}$ is zero for $v\in S$ and has degree $0$ for $v\notin S$ and
moreover property (3) of \cite[Prop. 7.5]{HW} holds (the relevance of this
property for us is the following: we can choose $\Phi$ as in loc. cit. so
that  property (3) will ensure that $z'_{\A}$ satisfies the assumptions (ii)
and(iii) of \cite[Thm. 6.2]{HW}). From our choice of $S$ it is clear that
$z'_{\A}$ is a family of effective, reduced local zero-cycles of constant
degree $\delta$ on $Y$, which is orthogonal to $B$.

We now define a Hilbert subset of $\P^1_k$ as follows.
According to \cite[Thm C]{OS} there exists a constant $C$ such that $$|\Br(\ov X)^{\Gal(\ov k/N)}| \leq C$$ for any field extension $N/k$ with $[N:k]\leq \delta$. Since $\Br(\ov X)\cong (\Q/\Z)^{22-\rho}$, where $\rho$ is the
geometric Picard number, there are finitely many elements of  $\Br(\ov X) $ whose order is less or equal to $C$. Each one of them is stabilized by a subgroup of $\Gamma_k$ of finite index, which corresponds to a finite extension
of $k$. Taking the Galois closure of the finitely many such extensions we produce a finite Galois extension $M/k$.
%Note that by construction if $K/k$ is an extension that is linearly disjoint from $M$ over $k$ then
%$A^{\Gamma_k}=A^{\Gamma_K}$, where $A$ is any $\Gamma_k$-invariant subgroup of $\Br(\ov X)$ consisting of elements of order less or equal to $C$.
We choose $L$ as in lemma \ref{theL} and denote by $F$ be the composite of $M$ and $L$. For this $F$ we now choose a Hilbert subset $H$ of $\P^1_k$ as in lemma \ref{lemhs}.

We are now in a position to apply \cite[Thm. 6.2]{HW}. Let $B$, $z'_{\A}$ and $H$ be as above. For these choices the
assumptions of loc. cit. are fullfiled and hence there exists a family
$\{z''_v\}_{v \in \Omega}\in Z_{0,\A}(Y)$ and a closed point $P\in H$ such
that $\{z''_v\}$ is effective for all $v$, $p_*(z''_v)=P\in Z_0(\P^1_{k_v})$
for all $v$, and $\{z''_v\}_{v \in \Omega}$ is orthogonal to $B$. In
particular $[k(P):k]=\delta$. Denote $k(P)$ by $k'$ and consider the fibre of $Y$ at $P$, i.e. $X_{k'}$. The image of $P$ in $Z_0(\P^1_{k_v})$ is the disjoint union $\sqcup_{w\in R} \Spec(k'_w)$ where $R$ is the set of the different prolongations of $v$ to $k'$. It is not difficult to see that because of the above properties of $z''_v$, we must have that $z''_v$ is the sum of distinct closed points, one closed point above each $\Spec(k'_w)$, and that closed point must have residue field $k'_w$. This implies that $\{z''_v\}_{v \in \Omega_k}$ induces a family of local points of the variety $X_{k'}$, i.e. an element of $\prod_{w\in \Omega_{k'}} X(k'_w)$. Furthermore, this family is orthogonal to the image of $B$ in
  $\Br(X_{k'})$.

 Assume that the image of $B$ under the natural map generates $\Br(X_{k'})/\Br_0(X_{k'})$. Since the Brauer-Manin obstruction to the Hasse principle on $X_{k'}$ is the only one, we deduce the existence of a
  $k'$-rational point of $X_{k'}$. This gives us a zero cycle, say $a$, of degree $\delta$ on $Y$. Then $a+z$ is a zero cycle of degree $d$ on $Y$.

  It remains to show that the image of $B$ generates $\Br(X_{k'})/\Br_0(X_{k'})$ or equivalently that the map $\Br(X)/\Br_0(X)\to\Br(X_{k'})/\Br_0(X_{k'})$ is surjective.
  Denote by $F'$ the composite of $F$ and $k'$. Note that $F'/k'$ is a Galois extension and there is a natural isomorphism $\Gal(F'/k') \to \Gal(F/k)$ compatible with the actions on $\Pic (\ov X)$.
  Moreover since $[k':k]=\delta$ we have that $|\Br( \ov X)^{\Gamma_{k'}}| \leq C$. This implies that $\Br( \ov X)^{\Gamma_{k}}=\Br( \ov X)^{\Gamma_{k'}}$: Let $\alpha \in \Br( \ov X)^{\Gamma_{k'}}$. The stabilizer of $\alpha$ in
  $\Gamma_k$ corresponds to a sub-extension of $F/k$ which is contained in $k'/k$. As $k'$ and $F$ are linearly disjoint over $k$ by construction, this sub-extension of $F/k$ must be $k/k$, i.e. $a$ is stabilized by $\Gamma_k$.
By our constructions we have the following commutative diagram with exact rows

\[
\begin{tikzcd}[column sep=tiny]
0\arrow{r}&H^1(\Gal(F/k), \Pic \ov X)\arrow{r}\arrow{d} &\Br(X)/\Br_0(X)\arrow{r}\arrow{d}
&\Br( \ov X)^{\Gamma_k}\arrow{d} \arrow{r}& H^2(\Gal(F/k), \Pic(\ov X))\arrow{d}\\
0\arrow{r}&H^1(\Gal(F'/k'), \Pic \ov X)\arrow{r}&\Br(X_{k'})/\Br_0(X_{k'})\arrow{r}&\Br( \ov X)^{\Gamma_{k'}}\arrow{r}&H^2(\Gal(F'/k'), \Pic(\ov X))
\end{tikzcd}
\]

 where the first, third and fourth vertical arrows are isomorphisms. Therefore the second vertical arrow is an isomorphism as well.
\qed

{\it Proof of Theorem \ref{WA}: } First we claim that for any $n\geq 1$ and any
finite subset $S\subset\Omega$, there exists a zero-cycle $b$ on $X$ of the
same degree such that $b=z_v$ in $CH_0(X_{k_v})/n $ for all $v\in S$:

Let $Y=X\times \P^1$ and let $p:Y\to X$ be the projection to the first
factor. As before (cf. lemma \ref{lempr}) it is clear that it suffices to
show the following: let $z_{\A}=\{z_v\}_{v \in \Omega}$ be a family of local
zero-cycles of degree $d$ on $Y$, which is orthogonal to $\Br(Y)$. Then for
any finite subset $S\subset\Omega$ and any $n\geq 0$ there exists a
zero-cycle $y$ on $Y$ of the same degree such that $y=z_v$ in
$CH_0(Y_{k_v})/n $ for all $v\in S$.

  We first enlarge $S$ as in the beginning of the proof of Theorem \ref{mainthm}. We now proceed exactly like the proof of Theorem \ref{mainthm}, keeping that notation. When we apply \cite[Thm. 6.2]{HW} we can choose the $z''_v$
  as close as we want to $z'_v$ in $Sym_{Y/k}(k_v)$ for $v\in S$ (this is part of \cite[Thm. 6.2]{HW} ).
  Applying weak approximation to the family induced by the $\{z''_v\}_{v \in \Omega_k}$ on $X_{k(P)}$, we can therefore take the zero cycle $a$ coming from a $k(P)$-rational point of $X_{k(P)}$ to be as close as we want to $z'_v$ in
  $Sym_{Y/k}(k_v)$ for $v\in S$. By \cite[Lemme 1.8]{Wit12} we can therefore ensure that $a=z'_v\in CH_0(Y_{k_v})/n$ for all $v\in S$. Take $b=a+z$. The claim is established.

  We now take  $S$ large enough so that for all $v\notin S$ we have that $X_{k_v}$ has a smooth projective model over $\O_v$ and moreover the characteristic of the residue field of $k_v$ does not divide $n$.
   We apply the claim to this choice of $S$. Hence there exists a zero-cycle $b$ on $X$ of the
same degree such that $b=z_v$ in $CH_0(X_{k_v})/n $ for all $v\in S$. If $v\notin S$, then since $b$ and $z_v$ have the same degree we deduce from Proposition \ref{Saito} that $b-z_v$ is divisible by $n$ in $CH_0(X_{k_v})$. Therefore $b=z_v$ in $CH_0(X_{k_v})/n $ for all $v\in \Omega$.

\qed

{\it Acknowledgment.} The author thanks Olivier Wittenberg for useful discussions on the contents of this note.

\end{document}